\documentclass[12pt]{amsart}
\usepackage[latin1]{inputenc}
\usepackage{amsxtra,amssymb,amsthm,amsmath,amscd,mathrsfs, epsfig, eufrak}
\usepackage[all]{xy}
\makeatletter
\@namedef{subjclassname@2010}{%
  \textup{2010} Mathematics Subject Classification}
\makeatother

\def \F {{\mathbb F}}

\def \Q {{\mathbb Q}}

\def \Z {{\mathbb Z}}

\font\cmm=cmmi10 at 14pt
\def\cmmp{\hbox{\cmm p}}
\font\cmmm=cmmi10 at 11pt
\def\cmmmp{\hbox{\cmmm p}}

\def \d {\,{\rm d}}

\def\GL{\hbox{{\rm GL}}}
\def\Li{\hbox{{\rm Li}}}

\def\le{\leqslant}
\def\ge{\geqslant}

\theoremstyle{plain}
\newtheorem{theorem}{Theorem}[section]
\newtheorem{lemma}[theorem]{Lemma}

\theoremstyle{remark}

\numberwithin{equation}{section}

\begin{document}

\title[The average exponent of elliptic curves modulo $\cmmmp$]
{The average exponent of elliptic curves modulo $\cmmp$}
\author{J. Wu}

\address{%School of Mathematics
%\\
%Shandong University
%\\
%Jinan
%\\
%Shandong 250100
%\\
%China
%\\
Institut Elie Cartan UMR 7502
\\
CNRS, Universit\'e de Lorraine, INRIA
\\
54506 Van\-d\oe uvre-l\`es-Nancy
\\
France}
\email{Jie.Wu@univ-lorraine.fr}

\date{\today}

\begin{abstract}
Let $E$ be an elliptic curve defined over $\Q$.
For a prime $p$ of good reduction for $E$, 
denote by $e_p$ the exponent of the reduction of $E$ modulo $p$.
Under GRH, we prove that there is a constant $C_E\in (0, 1)$ such that
$$
\frac{1}{\pi(x)} \sum_{p\le x} e_p
= \frac{1}{2} C_E x + O_E\big(x^{5/6} (\log x)^{4/3}\big)
$$
for all $x\ge 2$, where the implied constant depends on $E$ at most.
When $E$ has complex multiplication, 
the same asymptotic formula with a weaker error term
$O_E(1/(\log x)^{1/14})$ is established unconditionally.
These improve some recent results of Freiberg and Kurlberg.
\end{abstract}
\maketitle

%\addtocounter{footnote}{1}

\section{Introduction}

\smallskip

Let $E$ be an elliptic curve defined over $\Q$.
For a prime $p$ of good reduction for $E$
the reduction of $E$ modulo $p$ is an elliptic curve $E_p$ defined over the finite field $\F_p$
with $p$ elements.
Denote by $E_p(\F_p)$ the group of $\F_p$-rational points of $E_p$.
Its structure as a group, for example, the existence of large cyclic subgroups, especially of prime order, is of interest because of applications to elliptic curve cryptography \cite{Koblitz1987, Miller1986}.
It is well known that the finite abelian group $E_p(\F_p)$ has structure
\begin{equation}\label{Structure}
E_p(\F_p)\simeq (\Z/d_p\Z) \oplus (\Z/e_p\Z)
\end{equation}
for uniquely determined positive integers $d_p$ and $e_p$ with $d_p\mid e_p$.
Here $e_p$ is the size of the maximal cyclic subgroup of $E_p(\F_p)$, called the exponent of $E_p(\F_p)$.
The study about $e_p$ as a function of $p$ has received considerable attention 
\cite{Schoof1991, Duke2003, Co2003, CoMu2004}, where the following problems were considered:
\begin{itemize}
\item{lower bounds for the maximal values of $e_p$,}
\item{the frequency of $e_p$ taking its maximal value, 
i.e., the density of the primes $p$ for which $E_p(\F_p)$ is a cyclic group,}
\item{the smallest prime $p$ for which the group $E_p(\F_p)$ is cyclic (elliptic curve analogue of Linnik's problem).}
\end{itemize}
Very recently motivated by a question of Silverman, 
Freiberg and Kurlberg \cite{FK2012} investigated the average order of $e_p$.
Before stating their results, let us fixe some notation.
Given a positive integer $k$,
let $E[k]$ denote the group of $k$-torsion points of $E$ 
(called {\it the $k$-division group of $E$}) and 
let $L_k := \Q(E[k])$ be the field obtained by adjoining to $\Q$ the coordinates of the points of $E[k]$ 
(called {\it the $k$-division field of $E$}).
Write
\begin{equation}\label{defnLk}
n_{L_k} := [L_k : \Q].
\end{equation}
Denote by $\mu(n)$ the M\"obius function, by $\pi(x)$ the prime-counting function and
by $\zeta_{L_k}(s)$ the Dedekind zeta function associated with $L_k$, respectively.
Assuming the Generalized Riemann Hypothesis (GRH) for $\zeta_{L_k}(s)$ for all positive integers $k$,
Freiberg and Kurlberg \cite[Theorem 1.1]{FK2012} shew that
\begin{equation}\label{FK}
\frac{1}{\pi(x)} \sum_{p\le x} e_p
= \frac{1}{2} C_E x + O_E\big(x^{9/10} (\log x)^{11/5}\big)
\end{equation}
for all $x\ge 2$, 
where 
\begin{equation}\label{defCE}
C_E
:= \sum_{k=1}^{\infty} \frac{1}{n_{L_k}} \sum_{dm=k} \frac{\mu(d)}{m}
= \prod_p \bigg(1-\sum_{\nu=1}^{\infty} \frac{p-1}{p^\nu n_{L_{p^\nu}}}\bigg).
\end{equation}
The implied constant depends on $E$ at most.
When $E$ has complex multiplication (CM), they \cite[Theorem 1.2]{FK2012} also proved that \eqref{FK} holds unconditionally with a weaker error term
\begin{equation}\label{FKCM}
O_E\bigg(x\frac{\log_3x}{\log_2x}\bigg),
\end{equation}
where $\log_\ell$ denotes the $\ell$-fold iterated logarithm.

\vskip 2mm

The aim of this short note is to propose more precise result than \eqref{FK} and \eqref{FKCM}.

\begin{theorem}\label{thm}
Let $E$ be an elliptic curve over $\Q$.
\par
{\rm (a)}
Assuming GRH for the Dedekind zeta function $\zeta_{L_k}$ for all positive integers $k$, we have
\begin{equation}\label{thm(a)}
\frac{1}{\pi(x)} \sum_{p\le x} e_p
= \frac{1}{2} C_E x + O_E\big(x^{5/6} (\log x)^{4/3}\big).
\end{equation}

{\rm (b)}
If $E$ has CM, then we have unconditionally
\begin{equation}\label{thm(b)}
\frac{1}{\pi(x)} \sum_{p\le x} e_p
= \frac{1}{2} C_E x + O_E\bigg(\frac{x}{(\log x)^{1/14}}\bigg).
\end{equation}
Here $C_E$ is given as in \eqref{defCE} and the implied constants depend on $E$ at most.
\end{theorem}

{\bf Remark}.
(a)
Our proof of Theorem \ref{thm} is a refinement of Freiberg and Kurlberg's method \cite{FK2012}
with some simplification.

(b)
For comparison of \eqref{FK} and \eqref{thm(a)}, we have $\tfrac{9}{10}=0.9$ and $\tfrac{5}{6}=0.833\cdots$.

(c)
The quality of \eqref{thm(b)} can be compared with the following result of Kurlberg and Pomerance \cite[Theorem 1.2]{KP2012} concernng the multiplicative order of a number modulo $p$ :
Given a rational number $g\not=0, \pm 1$ and prime $p$ not dividing the numerator of $g$,
let $\ell_g(p)$ denote the multiplicative order of $g$ modulo $p$.
Assuming GRH for $\zeta_{\Q(g^{1/k}, {\rm e}^{2\pi{\rm i}/k})}(s)$ for all positive integers $k$, one has
$$
\frac{1}{\pi(x)} \sum_{p\le x} \ell_g(p)
= \frac{1}{2} C_g x + O\bigg(\frac{x}{(\log x)^{1/2-1/\log_3x}}\bigg),
$$
where $C_g$ is a positive constant depending on $g$.

\vskip 5mm

\section{Preliminary}

Let $E$ be an elliptic curve over $\Q$ with conductor $N_E$ and let $k\ge 1$ be an integer.
For $x\ge 1$, define
$$
\pi_E(x; k)
:= \sum_{\substack{p\le x\\ p\nmid N_E, \, k\mid d_p}} 1.
$$
The evaluation of this function will play a key role in the proof of Theorem \ref{thm}.
Using the Hasse inequality (see \eqref{Hasse} below),
it is not difficult to check that $p\nmid d_p$ for $p\nmid N_E$.
Thus the conditions $p\nmid N_E$ and $k\mid d_p$ are equivalent to $p\nmid kN_E$ and $k\mid d_p$,
that is $p\nmid kN_E$ and $E_p(\F_p)$ contains a subgroup isomorphic to $\Z/k\Z\times \Z/k\Z$.
Hence by \cite[Lemma 1]{Mu1983}, we have
$$
\sum_{\substack{p\le x\\ \text{$p$ splits completely in $L_k$}}} 1
= \pi_E(x; k) + O(\log(N_Ex)).
$$
In order to evaluate the sum on the left-hand side,
we need effective versions of the Chebotarev density theorem.
They were first derived by Lagarias and Odlyzko \cite{LaOd1979},
refined by Serre \cite{Serre1981}, and subsequently improved by M. Murty, V. Murty and Saradha \cite{MuMuSa1988}.
With the help of these results, one can deduce the following lemma 
(cf. \cite[Lemma 3.3]{FK2012}).

\begin{lemma}\label{lem1}
Let $E$ be an elliptic curve over $\Q$ with conductor $N_E$.
\par
{\rm (a)}
Assuming GRH for the Dedekind zeta function $\zeta_{L_k}(s)$, we have
\begin{equation}\label{Eq1lem1}
\pi_E(x; k) = \frac{\Li(x)}{n_{L_k}} + O\big(x^{1/2} \log(N_Ex)\big)
\end{equation}
uniformly for $x\ge 2$ and $k\ge 1$,
where the implied constant is absolute.
\par
{\rm (b)}
There exist two absolute constants $B>0$ and $C>0$ such that
\begin{equation}\label{Eq2lem1}
\pi_E(x; k)= \frac{\Li(x)}{n_{L_k}} + O\big(x {\rm e}^{-B(\log x)^{5/14}}\big)
\end{equation}
unformly for $x\ge 2$ and $C N_E^2 k^{14}\le \log x$,
where the implied constant is absolute.
\end{lemma}

The next lemma (cf. \cite[Proposition 3.2]{FK2012} or \cite[Propositions 3.5 and 3.6]{CoMu2004})
gathers some properties of the division fields $L_k$ of $E$ and estimates for $n_{L_k}$,
which will be useful later. 
Denote by $\varphi(k)$ the Euler function.

\begin{lemma}\label{lem2}
{\rm (a)}
The field $L_k$ contains $\Q({\rm e}^{2\pi {\rm i}/k})$.
Therefore $\varphi(k)\mid n_{L_k}$ and
a rational prime $p$ which splits completely in $L_k$ satisfies $p\equiv 1 ({\rm mod}\,k)$.

\par

{\rm (b)}
$n_{L_k}$ divides $|\GL_2(\Z/k\Z)|=k^3 \varphi(k) \prod_{p\mid k} (1-p^{-2})$.

\par

{\rm (c)}
If $E$ is a non-CM curve, then there exists a constant $B_E\ge 1$ (depending only on $E$)
such that $|\GL_2(\Z/k\Z)|\le B_E n_{L_k}$ for each $k\ge 1$.
Moreover, we have $|\GL_2(\Z/k\Z)| = n_{L_k}$ whenover $(k, M_E)=1$ 
$($where $M_E$ is Serre's constant\,$)$.

\par
{\rm (d)}
If $E$ has CM, then $\varphi(k)^2\ll n_{L_k}\le k^2$.
\end{lemma}

\vskip 5mm

\section{Proof of Theorem \ref{thm}}

Let $a_E(p) := p+1- |E_p(\F_p)|$, then
$$
e_p = \begin{cases}
(p+1-a_E(p))/d_p & \text{if $\,p\nmid N_E$},
\\\noalign{\vskip 1mm}
0 & \text{otherwise}.
\end{cases}
$$
By using Hasse's inequality
\begin{equation}\label{Hasse}
|a_E(p)|<2\sqrt{p}
\end{equation}
for all primes $p\nmid N_E$, it is easy to see that
\begin{equation}\label{Eq1}
\sum_{p\le x} e_p
= \sum_{p\le x, \, p\nmid N_E} \frac{p}{d_p} + O\bigg(\frac{x^{3/2}}{\log x}\bigg).
\end{equation}

In order to evaluate the last sum,
we first notice that the Hasse inequality \eqref{Hasse} implies $d_p\le 2\sqrt{p}$.
Thus we can use the formula 
$$
\frac{1}{k} 
= \sum_{dm\mid k} \frac{\mu(d)}{m}
$$
to write
\begin{equation}\label{Eq2}
\sum_{\substack{p\le x\\ p\nmid N_E}} \frac{p}{d_p} 
= \sum_{\substack{p\le x\\ p\nmid N_E}} p \sum_{dm\mid d_p} \frac{\mu(d)}{m}
= \sum_{k\le 2\sqrt{x}} \sum_{dm=k} \frac{\mu(d)}{m} 
\sum_{\substack{p\le x\\ p\nmid N_E, \, k\mid d_p}} p.
\end{equation}
Let $y\le 2\sqrt{x}$ be a parameter to be choosen later and define
\begin{align*}
S_1
& := \sum_{k\le y} \sum_{dm=k} \frac{\mu(d)}{m} \sum_{\substack{p\le x\\ p\nmid N_E, \, k\mid d_p}} p,
\\
S_2
& := \sum_{y<k\le 2\sqrt{x}} \sum_{dm=k} \frac{\mu(d)}{m} \sum_{\substack{p\le x\\ p\nmid N_E, \, k\mid d_p}} p.\end{align*}
With the help of Lemma \ref{lem1}(a), a simple partial integration allows us to deduce 
(under GRH)
\begin{equation}\label{sump}
\begin{aligned}
\sum_{\substack{p\le x\\ p\nmid N_E, \, k\mid d_p}} p
& = \int_{2-}^x t \d \pi_E(t; k)
= x \pi_E(x; k)
- \int_{2}^x \pi_E(t; k) \d t 
\\\noalign{\vskip -3mm}
& = \frac{x \Li(x)}{n_{L_k}} 
- \frac{1}{n_{L_k}} \int_{2}^x \Li(t) \d t 
+ O_E\big(x^{3/2}\log x\big)
\\\noalign{\vskip 1mm}
& = \frac{\Li(x^2)}{n_{L_k}} 
+ O_E\big(x^{3/2}\log x\big).
\end{aligned}
\end{equation}
On the other hand, by Lemma \ref{lem2} we infer that
\begin{equation}\label{CE}
\sum_{k\le y} \frac{1}{n_{L_k}} \sum_{dm=k} \frac{\mu(d)}{m}
= C_E + O(y^{-1}).
\end{equation}
Thus combining \eqref{sump} with \eqref{CE} and using the following trivial inequality
\begin{equation}\label{trivial}
\bigg|\sum_{dm=k} \frac{\mu(d)}{m}\bigg|
\le \frac{\varphi(k)}{k}
\le 1,
\end{equation}
we find
\begin{equation}\label{S1}
\begin{aligned}
S_1
& = \Li(x^2) \sum_{k\le y} \frac{1}{n_{L_k}}\sum_{dm=k} \frac{\mu(d)}{m}
+ O_E\bigg(x^{3/2}\log x \sum_{k\le y} \bigg|\sum_{dm=k} \frac{\mu(d)}{m}\bigg|\bigg)
\\
& = C_E \Li(x^2) 
+ O_E\bigg(\frac{x^2}{y\log x} + x^{3/2}y\log x\bigg).
\end{aligned}
\end{equation}

Next we treat $S_2$.
By \cite[Lemma 3.1 and Proposition 3.2(a)]{FK2012},
we see that $k\mid d_p$ implies that $k^2\mid (p+1-a_E(p))$ and also $k\mid (p-1)$,
hence $k\mid (a_E(p)-2)$.
With the aid of this and the Brun-Titchmarsh inequality, we can deduce that
\begin{align*}
S_2
& \ll x \sum_{y<k\le 2\sqrt{x}} 
\bigg(
\sum_{\substack{|a|\le 2\sqrt{x}, a\not=2\\ a\equiv 2 ({\rm mod} k)}}
\sum_{\substack{p\le x, a_E(p)=a\\ k^2\mid p+1-a}} 1
+
\sum_{\substack{p\le x, a_E(p)=2\\ k^2\mid p-1}} 1
\bigg)
\\
& \ll x \sum_{y<k\le 2\sqrt{x}} 
\bigg(
\frac{\sqrt{x}}{k}\cdot\frac{x}{k\varphi(k)\log(8x/k^2)}
+
\frac{x}{k^2}
\bigg).
\end{align*}
By virtue of the elementary estimate
$$
\sum_{n\le t} \frac{1}{\varphi(k)}
= D\log t + O(1)
\qquad(t\ge 1)
$$
with some positive constant $D$,
a simple integration by parts leads to
\begin{equation}\label{S2}
S_2
\ll \frac{x^{5/2}}{y^2 \log(8x/y^2)} + \frac{x^{2}}{y}\cdot
\end{equation}
Inserting \eqref{S1} and \eqref{S2} into \eqref{Eq2}, we find
\begin{equation}\label{Eq3}
\begin{aligned}
\sum_{p\le x, \, p\nmid N_E} \frac{p}{d_p} 
= C_E \Li(x^2) 
+ O_E\bigg(x^{3/2}y\log x + \frac{x^{5/2}}{y^2 \log(8x/y^2)} + \frac{x^{2}}{y}\bigg),
\end{aligned}
\end{equation}
where we have used the fact that the term $x^2y^{-1}(\log x)^{-1}$ 
can be absorded by $x^{5/2}y^{-2}(\log(8x/y^2))^{-1}$
since $y\le 2\sqrt{x}$.
Now the asymptotic formula \eqref{thm(a)} follows from \eqref{Eq1} and \eqref{Eq3} 
with the choice of $y=x^{1/3}(\log x)^{-2/3}$.

\vskip 1mm

The proof of \eqref{thm(b)} is very similar to that of \eqref{thm(a)}.
Next we shall only point out some important differences.

Similar to \eqref{sump}, we can apply Lemma \ref{lem1}(b) to prove (unconditionally)
$$
\sum_{\substack{p\le x\\ p\nmid N_E, \, k\mid d_p}} p
= \frac{\Li(x^2)}{n_{L_k}} 
+ O_E\big(x^{2}\exp\{-B(\log x)^{5/14}\}\big)
$$
for $k\le (C^{-1}N_E^{-2}\log x)^{1/14}$.
As before from this and \eqref{CE}-\eqref{trivial}, we can deduce that
\begin{equation}\label{S1(b)}
S_1
= C_E \Li(x^2) 
+ O_E\big(x^2y^{-1}(\log x)^{-1} + x^{2}y{\rm e}^{-B(\log x)^{5/14}}\big)
\end{equation}
for $y\le (C^{-1}N_E^{-2}\log x)^{1/14}$.

The treatment of $S_2$ is different.
First we divide the sum over $k$ in $S_2$ into two parts accroding to
$y<k\le x^{1/4}(\log x)^{3/4}$ or $x^{1/4}(\log x)^{3/4}<k\le 2\sqrt{x}$.

When $E$ has CM, we have (see \cite[page 692]{Duke2003})
$$
\sum_{\substack{p\le x\\ p\nmid N_E, \, k\mid d_p}} 1
\ll \frac{x}{\varphi(k)^2 \log x}
$$
for $k\le x^{1/4}(\log x)^{3/4}$.
Thus the contribution from $y<k\le x^{1/4}(\log x)^{3/4}$ to $S_2$ is
\begin{align*}
& \ll \frac{x^2}{\log x} \sum_{y<k\le x^{1/4}(\log x)^{3/4}} \frac{1}{\varphi(k)^2}
\ll \frac{x^2}{y\log x}\cdot
\end{align*}
Clearly the inequality \eqref{S2} (taking $y=x^{1/4}(\log x)^{3/4}$) implies that the contribution from 
$x^{1/4}(\log x)^{3/4}<k\le 2\sqrt{x}$ to $S_2$ is
\begin{align*}
\ll \sum_{x^{1/4}(\log x)^{3/4}<k\le 2\sqrt{x}}
\sum_{\substack{p\le x\\ p\nmid N_E, \, k\mid d_p}} p
& \ll \frac{x^2}{(\log x)^{5/2}}\cdot
\end{align*}
By combining these two estimates, we obtain
\begin{equation}\label{S2(b)}
S_2\ll \frac{x^2}{y\log x} + \frac{x^2}{(\log x)^{5/2}}\cdot
\end{equation}
Inserting \eqref{S1(b)} and \eqref{S2(b)} into \eqref{Eq2}, we find
\begin{equation}\label{Eq3(b)}
\begin{aligned}
\sum_{p\le x, \, p\nmid N_E} \frac{p}{d_p} 
= C_E \Li(x^2) 
+ O_E\bigg(\frac{x^2}{y\log x} + \frac{x^2}{(\log x)^{5/2}} + x^{2}y{\rm e}^{-B(\log x)^{5/14}}\bigg)
\end{aligned}
\end{equation}
for $y\le (C^{-1}N_E^{-2}\log x)^{1/14}$.

Now the asymptotic formula \eqref{thm(b)} follows from \eqref{Eq1} and \eqref{Eq3(b)} 
with the choice of $y = (C^{-1}N_E^{-2}\log x)^{1/14}$.


\vskip 10mm

\begin{thebibliography}{CC}

\bibitem{Co2003}
A. C. Cojocaru,
\textit{Cyclicity of CM elliptic curves modulo $p$},
Trans. AMS. {\bf 355} (2003), 2651--2662.

\bibitem{CoMu2004}
A. C. Cojocaru and M. Ram Murty,
\textit{Cyclicity of elliptic curves modulo $p$ and elliptic curve analogue of Linnik's problem},
Math. Ann. {\bf 330} (2004), no. 7, 601--625.

\bibitem{Duke2003}
W. Duke,
\textit{Almost all reductions modulo $p$ of an elliptic curve have a large exponent}, 
C. R. Acad. Sci. Paris, Ser. {\bf 1337} (2003), 689--692.

\bibitem{FK2012}
T. Freiberg and P. Kurlberg,
\textit{On the average exponent of elliptic curves modulo $p$}, 
arXiv:1203.4382v1. (21 pages)

\bibitem{Koblitz1987}
N. Koblitz,
\textit{Elliptic curve cryptosystems}, 
Math. Comp. {\bf 48} (1987), 203--209.

\bibitem{KP2012}
P. Kurlberg and C. Pomerance,
\textit{On a problem of Arnold: the average multiplicative order of a given integer},
to appear in Algebra and Number Theory.

\bibitem{LaOd1979}
J. Lagarias and A. Odlyzko,
\textit{Effective versions of the Chebotarev Density Theorem},
in: Algebraic Number Fields (A. Fr\"ohlich edit.), New York, Academic Press (1977), 409--464.

\bibitem{Miller1986}
V. S. Miller,
\textit{Use of elliptic curves in cryptography}, 
Advances in cryptology --- CRYPTO 85 (Santa Barbara, Calif., 1885),
Lecture Notes in Comput. Sci. {\bf 218}, Springer, Berlin, 1986, 417--426.

\bibitem{Mu1983}
M.-R. Murty,
\textit{On Artin's conjecture},
J. Number Theory {\bf 16} (1983), 147--168.

\bibitem{MuMuSa1988}
M.-R. Murty, V.-K. Murty and N. Saradha,
\textit{Modular forms and the Chebotarev density theorem},
Amer. J. Math. {\bf 110} (1988), 253--281.

\bibitem{Schoof1991}
R. Schoof,
\textit{The exponent of the group of points on the reductions of an elliptic curve}, 
in : Arithmetic algebraic geometry (Texel, 1989),
in : Progr. Math. {\bf 89} (1991), Birkh\"auser, Boston, MA, 325--335.

\bibitem{Serre1981}
J.-P. Serre,
\textit{Quelques applications du th\'eor\`eme de densit\'e de Chebotarev},
Inst. Hautes Etudes Sci. Publ. Math. {\bf 54} (1981), 123--201.


\end{thebibliography}
\end{document}